\documentclass[12pt]{amsart}

\usepackage{mathrsfs}

\usepackage{amsmath,amsthm,amssymb}
\font\teneufm=eufm10
\font\seveneufm=eufm7
\font\fiveeufm=eufm5
\newfam\frakturfam

\textfont\frakturfam=\teneufm
\scriptfont\frakturfam=\seveneufm
\scriptscriptfont\frakturfam=\fiveeufm

\newtheorem{pr}{Proposition}

\newtheorem{lm}{Lemma}
\newtheorem{theor}{Theorem}
\newtheorem{co}{Corollary}

\def\bee{\begin{eqnarray}}
\def\bes{\begin{eqnarray*}}
\def\eee{\end{eqnarray}}
\def\ees{\end{eqnarray*}}

\def\a{\alpha}
\def\b{\beta}

\def\Proof{{\sl Proof.}\ }

\title{Algorithmic problems for differential polynomial algebras}

\begin{document}
\date{}
\maketitle

\begin{center}
{\bf Ualbai Umirbaev}\footnote{Supported by an MES grant 1226/GF3; Eurasian National University,
 Astana, Kazakhstan and
 Wayne State University,
Detroit, MI 48202, USA,
e-mail: {\em umirbaev@math.wayne.edu}}
\end{center}

\begin{abstract} We prove that the ideal membership problem and the subalgebra membership problem are algorithmically undecidable for differential polynomial algebras with at least two basic derivation operators.
\end{abstract}

\noindent


\section{Introduction}

\hspace*{\parindent}

Let $P_n=k[x_1,x_2,\ldots,x_n]$ be the polynomial algebra in the variables $x_1,x_2,\ldots,x_n$ over a constructive field $k$. One of the first applications of Gr\"obner bases (see, for example \cite{CLO}) gives the decidability of the ideal membership problem for $P_n$, i.e., there exists an effective algorithm which for any finite sequence of elements $f,f_1,\ldots,f_m\in P_n$ determines whether $f$ belongs to the ideal $(f_1,\ldots,f_m)$ or not. Another application of Gr\"obner bases \cite{SS} gives the decidability of the subalgebra membership problem for $P_n$, i.e., there exists an effective algorithm which for any finite sequence of elements $f,f_1,\ldots,f_m\in P_n$ determines whether $f$ belongs to the subalgebra $\langle f_1,\ldots,f_m\rangle$ or not. The subalgebra membership problem in characteristic zero was also solved in \cite{Noskov} without Gr\"obner bases.

Traditionally, the ideal membership problem for free algebras is called the word problem for corresponding variety of algebras. The word problem is undecidable for many subvarities of semigroups \cite{Gurevich}, groups, and associative and Lie algebras \cite{Kharlampovich81,SK}. More details on this classical problem can be found in a survey paper \cite{KS}. The decidability of the word problem, in general, is related to the study of Gr\"obner-Shirshov bases \cite{BC}. The word problem is decidable for polynilpotent $\mathfrak{N}_2\mathfrak{A}$-groups \cite{Kharlampovich87} and for polynilpotent $\mathfrak{N}_2\mathfrak{N}_c$-Lie algebras \cite{BU}.

A well known Nielsen-Schreier Theorem states that the subgroups of free groups are free \cite{KMS} and a Shirshov-Witt Theorem states that the subalgebras of free Lie algebras are free \cite{Shir1,Witt}. These results easily imply the decidability of the subalgebra membership problem for free groups and free Lie algebras. The subalgebra membership problem is decidable also for free metabelian groups \cite{Romanovskii} and free metabelian Lie algebras \cite{Zaicev}. It is undecidable for free associative algebras \cite{Um7} and for free solvable Lie algebras \cite{Um8} and for free solvable groups \cite{Um14} of solvability index $\geq 3$. The subalgebra membership problem for free metanilpotent Lie algebras, i.e.,  $\mathfrak{N}_s\mathfrak{N}_t$-Lie algebras,  is decidable \cite{GU1,GU2}.

The basic concepts of differential algebras can be found in \cite{Kolchin,Ritt,PS03}.
Let $\Delta = \{\delta_1,\ldots, \delta_m\}$ be a basic set of derivation operators and
let $\Phi\{x_1,x_2,\ldots,x_n\}$ be the differential polynomial ring in free differential variables $x_1,x_2,\ldots,x_n$ over an arbitrary commutative ring $\Phi$ with unity such that $\delta_i(\Phi)=0$ for all $i$.
There are several approaches to define analogues of the Gr\"obner bases for differential polynomial algebras \cite{CF89,Man91,Ol91} and some recent results can be found in \cite{KLMP}. The differential ideal membership problem is solved positively only in some interesting particular cases (see, for example in \cite{KZ,Z}). At the moment the membership problem for differential ideals generated by a single polynomial is still open. It is negatively solved for recursively generated differential ideals \cite{GMO}. The membership problem for finitely generated differential ideals of differential polynomial algebras was formulated by J.F. Ritt in \cite[p. 177, Question 2]{Ritt}.

In this paper we prove that the membership problem for finitely generated differential ideals is algorithmically undecidable, i.e., the word problem for differential algebras is undecidable. The main instrument of proving this is an interpretation of Minsky machines. The proof uses the fact that every recursive function can be calculated by Minsky machines without cycles \cite{Um7}. Using a method of interpreting  the ideal membership problem from \cite{U95A}, we also prove that the membership problem for finitely generated differential subalgebras is undecidable.

Our proofs need at least two derivation operators. Thus, these problems are still open for ordinary differential polynomial algebras.

The rest of the paper is organized as follows. In Section 2 we fix some standard notations and recall some definitions on differential algebras. In Section 3, we give an interpretation of the Minsky machines and prove the undecidability of the ideal membership problem. In Section 4 we give an interpretation of the ideal membership problem and prove  the undecidability of the subalgebra membership problem.

\section{Differential polynomial algebras}

\hspace*{\parindent}

All our rings are assumed to be commutative and with unity. Let $R$ be an arbitrary ring. A mapping $d: R\rightarrow R$ is called a {\em derivation} if
\bes
d(s+t)=d(s)+d(t)
\ees
\bes
d(st)=d(s)t+sd(t)
\ees
holds for all $s,t\in R$.

Let $\Delta = \{\delta_1,\ldots, \delta_m\}$ be a basic set of derivation operators.

A ring $R$ is said to be a {\em differential} ring or $\Delta$-ring if all elements of $\Delta$ act on $R$ as a commuting set of derivations, i.e., the derivations $\delta_i: R\rightarrow R$ are defined for all $i$ and $\delta_i\delta_j=\delta_j\delta_i$ for all $i,j$.

Let $\Theta$ be the free commutative monoid on the set $\Delta = \{\delta_1,\ldots, \delta_m\}$  of
derivation operators. The elements
\bes
\theta = \delta_1^{i_1}\ldots \delta_m^{i_m}
\ees
of the monoid $\Theta$ are called {\em derivative}
operators. The {\em order} of $\theta$ is defined as $|\theta|=i_1+\ldots+i_m$.

Let $R$ be a differential ring. Denote by $R^e$ the free left $R$-module with a basis $\Theta$. Every element $u\in R^e$ can be uniquely written in the form
\bes
u=\sum_{\theta\in \Theta} r_{\theta}\theta
\ees
with a finite number of nonzero $r_{\theta}\in R$. We turn $R^e$ to a ring by
\bes
\delta_i r=r \delta_i+ \delta_i(r)
\ees
for all $i$ and $r\in R$. It is well known \cite{KLMP} that these relations uniquely define a structure of a ring on $R^e$ and every left module over $R^e$ is called a {\em differential module} over $R$. In particular, $R$ is a left $R^e$ and every $I\subseteq R$ is a differential ideal of $R$ if and only if $I$ is an $R^e$-submodule of $R$. The ring $R^e$ is called the {\em universal enveloping} ring of $R$.

Let $x^\Theta=\{x^\theta | \theta\in \Theta\}$ be a set of symbols enumerated by the elements of $\Theta$. Consider the polynomial algebra $R[x^\Theta]$ over $R$ generated by the set of (polynomially) independent variables $x^\Theta$. It is easy to check that the derivations $\delta_i$ can be uniquely extended to a derivation of $R[x^\Theta]$ by $\delta_i(x^\theta)=x^{\delta_i\theta}$. Denote this differential ring
by $R\{x\}$; it is called the {\em ring of differential polynomials} in $x$ over $R$.

By adjoining more variables, we can obtain the differential ring $R\{x_1,x_2,\ldots,x_n\}$ of the differential polynomials in $x_1,x_2,\ldots,x_n$ over $R$. Let $M$ be the free commutative monoid generated by all elements $x_i^{\theta}$, where $1\leq i\leq n$ and $\theta\in \Theta$. The elements of $M$ are called {\em monomials} of $R\{x_1,x_2,\ldots,x_n\}$. Every element
$a\in R\{x_1,x_2,\ldots,x_n\}$ can be uniquely written in the form
\bes
a=\sum_{m\in M} r_m m
\ees
with a finite number of nonzero $r_m\in R$.

Every ring can be considered as a differential ring under the trivial action of all derivation operators. If all differential operators act as zeroes on $R$, then  $R\{x_1,x_2,\ldots,x_n\}$ becomes an $R$-algebra. In studying of Gr\"obner bases, we usually assume that $R$ is a constructive field $k$ or the ring of integers $\mathbb{Z}$.

\section{The ideal membership problem}

\hspace*{\parindent}

Minsky machines are multi-tape Turing machines \cite{Malcev}.
The hardware of a two-tape Minsky machine consists of two tapes and a head. The tapes are infinite to the right and are divided into infinitely many cells numbered from the left to the right, starting with zero.  The external alphabet consists of $0$ and $1$. The first cells on both tapes always contain $1$ and all other cells have $0$. The head may acquire one of several internal states: $q_0,q_1,\ldots,q_n$; $q_0$ is the {\em terminal} state. At every moment the head looks at one cell of the first tape and at one cell of the second tape.

The program of a Minsky machine consists of a set of commands of the form
\bee\label{f1}
q_i \varepsilon \sigma \rightarrow q_j T_{\a}T_{\b},
\eee
where $1\leq i\leq n$, $0\leq j\leq n$, $\varepsilon,\sigma\in \{0,1\}$, $\a,\b\in \{-1,0,1\}$, and $\a\geq 0$ if $\varepsilon=1$ and $\b\geq 0$ if $\sigma=1$. This means that if the head is in the state $q_i$ observing a cell containing $\varepsilon$ on the first tape and a cell containing $\sigma$ on the second tape, then it acquires the state $q_j$ and the first (the second) tape is shifted $\a$ (resp. $\b$) cells to the left relative to the head. If $\a=-1$, for example, then the first tape is shifted one cell to the right.

 A configuration of a Minsky machine can be described by a triple $[i,m,n]$, where $m$ and $n$ are the numbers of the cells observed by the head in the first and the second tapes, respectively, and $q_i$ is the internal state of the head. We write
\bes
[i,m,n]\rightarrow [j,p,q],
\ees
if a Minsky machine at the configuration $[i,m,n]$ gets the configuration [j,p,q] in one step, i.e., as a result of execution of one (a unique!) command of the type (\ref{f1}).

Recall that in algorithmic theory the set of natural numbers includes $0$, i.e., $\mathbb{N}=\{0,1,2,\ldots\}$.
Minsky \cite{Malcev} proved that for every partial recursive function $f: \mathbb{N}\rightarrow \mathbb{N}$ there exists a Minsky machine that calculates $f(x)$, i.e., for every natural $x$ it passes from the configuration $[1,2^x,0]$ to the configuration $[0,2^{f(x)},0]$ if $f(x)$ is defined, and operates infinitely, never reaching the terminal state $q_0$, if $f(x)$ is not defined.

We say that a Minsky machine has a cycle if there exists a configuration $[i,m,n]$ such that the machine starting work at this configuration returns to the same configuration in a finite number of positive steps. A Minsky machine without cycles is called {\em acyclic}.

We need the next lemma.
\begin{lm}\label{l1} \cite{Um7}
 Let $S$ be a recursively enumerable subset of natural numbers $\mathbb{N}$. Then there exists a two-tape acyclic Minsky machine that for every $x\in \mathbb{N}$ starting work at the configuration $[1,2^{2^x},0]$ reaches $[0,1,0]$ in finitely many steps if $x\in S$ and operates infinitely if $x\notin S$.
\end{lm}

First of all we assume that the basic set of derivations $\Delta = \{\delta_1,\ldots, \delta_m\}$ contains at least two elements. Moreover, we may assume that $\Delta = \{\delta_1,\delta_2\}$ since the other derivations do not hurt our proofs.

We also fix a recursively enumerable subset  $S$  of the set of natural numbers $\mathbb{N}$ and fix an acyclic  Minsky machine $M$ that calculates the characteristic function of $S$ as in Lemma \ref{l1}. Assume that (\ref{f1}) is the set of all commands of $M$.

Let $\Phi$ be an arbitrary ring. We consider all our algebras over $\Phi$. In the case of positive solutions of algorithmic problems we have to assume that $\Phi$ is constructive (or computable). But it is not mandatory for negative solutions. Of course, we assume that $\Phi$ contains a nonzero unity.

We consider $\Phi$ as a differential ring with the trivial action of the derivation operators. Let $A=\Phi\{x_1,x_2,q_0,q_1,\ldots,q_n\}$ be the free differential algebra over $\Phi$ in free differential variables $x_1,x_2,q_0,q_1,\ldots,q_n$.

With each command of $M$ of the type (\ref{f1}) we associate the element
\bes
f(i,\varepsilon,\sigma)=x_1^{\varepsilon}x_2^{\sigma}\delta_1^{1-\varepsilon}
\delta_2^{1-\sigma}(q_i)
-x_1^{\varepsilon}x_2^{\sigma}\delta_1^{1-\varepsilon+\a}
\delta_2^{1-\sigma+\b}(q_j)
\ees
of the algebra $A$, where $1\leq i\leq n$ and $\varepsilon, \sigma= 0,1$. Denote by $I$ the differential ideal of $A$ generated by all elements $f(i,\varepsilon,\sigma)$.

Denote by $J$ the differential ideal of $A$ generated by the elements
\bes
\delta_1(x_2), \delta_2(x_1).
\ees

Put also
\bes
f_m=x_1x_2\delta_1^{2^{2^m}}(q_1)-x_1x_2\delta_1(q_0)
\ees
for all $m\in \mathbb{N}$.

\begin{pr}\label{p1}
Element $f_m$ of $A$ belongs to the differential ideal $I+J$ if and only if $m\in S$.
\end{pr}

The rest of this section is devoted to the proof of this proposition.

Denote by $B$ the quotient algebra $A/J$.
\begin{lm}\label{l2} The algebra $B$ is a polynomial algebra over $\Phi$ in the polynomial variables
\bee\label{f2}
\delta_1^i(x_1), 
\delta_2^i(x_2), \delta_1^i\delta_2^j(q_0),\ldots, \delta_1^i\delta_2^j(q_n),
\eee
where $i,j\geq 0$.
\end{lm}
\Proof Let $R$ be a polynomial algebra over $\Phi$ in the set of variables (\ref{f2}). We can turn $R$ to a differential algebra by
\bes
\delta_1(\delta_1^i(x_1))=\delta_1^{i+1}(x_1), \delta_2(\delta_1^i(x_1))=0, 
\delta_2(\delta_2^i(x_2))=\delta_1^{i+1}(x_2),\\ \delta_1(\delta_2^i(x_2))=0, 
\delta_1(\delta_1^i\delta_2^j(y))=\delta_1^{i+1}\delta_2^j(y), 
\delta_2(\delta_1^i\delta_2^j(y))=\delta_1^i\delta_2^{j+1}(y), 
\ees
for all $i,j\geq 0$ and $y\in \{q_0,\ldots,q_n\}$.

Consider the differential homomorphism $\varphi: A\rightarrow R$ defined by $\varphi(x)=x$ for all $x=x_i,q_i$. Obviously, $\varphi(J)=0$ and it easy to check that the induced homomorphism $A/J\rightarrow R$ is an isomorphism. $\Box$

We continue to work with the algebra $B$. The images of elements of $A$ in $B$ will be written in the same way as in the algebra $A$. The images of $f(i,\varepsilon,\sigma), f_m$, and $I$ will be denoted by  $g(i,\varepsilon,\sigma), g_m$, and $\widetilde{I}$, respectively.

Notice that $B$ is homogeneous with respect to each of its polynomial generators (\ref{f2}). Moreover, the elements $g(i,\varepsilon,\sigma)$ and $g_m$ are homogeneous with respect to
each of the polynomial variables
\bee\label{f3}
\delta_1^i(x_1), \delta_2^i(x_2), \ \ i\geq 0,
\eee
and with respect to the group of variables
\bee\label{f5}
\delta_1^i\delta_2^j(q_0),\ldots, \delta_1^i\delta_2^j(q_n), \ \ i,j\geq 0.
\eee

Denote by $V$ the set of all monomials in the set of commuting variables (\ref{f2}). Every element of the universal enveloping algebra $B^e$ can be uniquely represented as a linear
combination of elements of the form
\bee\label{f6}
v\delta_1^i\delta_2^j, \ \ v\in V, i,j\geq 0.
\eee

Let $\deg$ be the standard polynomial degree function on $B$, i.e., $\deg(y)=1$ for all elements from (\ref{f2}). All elements $g(i,\varepsilon,\sigma)$ and $g_m$ are homogeneous with respect to $\deg$ and
\bes
\deg(g(i,\varepsilon,\sigma))=1+\varepsilon+\sigma, \ \ \deg(g_m)=3.
\ees

We also define polynomial degree functions $\deg_1$ and $\deg_2$ on $B$ as follows:
$\deg_1(\delta_1^i(x_1))=i+1$ for all $i\geq 0$ and $\deg(y)=0$ for all other variables from (\ref{f2}); $\deg_2(\delta_2^j(x_2))=j+1$ for all $j\geq 0$ and $\deg(y)=0$ for all other variables from (\ref{f2}). For any $v\in V$ put $\mathrm{Deg}(v)=(\deg_1(v),\deg_2(v))$. Let $\leq$ be the lexicographic order on $\mathbb{N}^2$ (recall that $\mathbb{N}$ includes $0$).
For any $v\in V$ denote by $\overline{v}$ its highest homogeneous part with respect to $\mathrm{Deg}$.
 The elements $g(i,\varepsilon,\sigma)$ and $g_m$ are also homogeneous with respect to  $\mathrm{Deg}$.

\begin{lm}\label{l3}
\bes
\overline{\delta_1^s\delta_2^tg(i,\varepsilon,\sigma)}
=(\delta_1^s(x_1))^{\varepsilon}(\delta_2^t(x_2))^{\sigma}\delta_1^{(s+1)(1-\varepsilon)}
\delta_2^{(t+1)(1-\sigma)}(q_i)\\
-(\delta_1^s(x_1))^{\varepsilon}(\delta_2^t(x_2))^{\sigma}\delta_1^{(s+1)(1-\varepsilon)+\a}
\delta_2^{(t+1)(1-\sigma)+\b}(q_j).
\ees
\end{lm}
\Proof We consider only the case $\varepsilon=1$ and $\sigma=0$ since the other cases can be treated similarly. We have
\bes
g(i,1,0)=x_1\delta_2^1(q_i)
-x_1\delta_1^{\a}\delta_2^{1+\b}(q_j).
\ees
Consequently,
\bes
\delta_1^s\delta_2^tg(i,1,0)=\delta_1^s\delta_2^t(x_1\delta_2^1(q_i)
-x_1\delta_1^{\a}\delta_2^{1+\b}(q_j))\\
=\delta_1^s(x_1\delta_2^{t+1}(q_i)
-x_1\delta_1^{\a}\delta_2^{t+1+\b}(q_j))\\
=\sum_{r=0}^s(\delta_1^r(x_1)\delta_1^{s-r}\delta_2^{t+1}(q_i)
-\delta_1^r(x_1)\delta_1^{s-r+\a}\delta_2^{t+1+\b}(q_j)).
\ees
Consequently,
\bes
\overline{\delta_1^s\delta_2^tg(i,1,0)}=\delta_1^s(x_1)\delta_2^{t+1}(q_i)
-\delta_1^s(x_1)\delta_1^{\a}\delta_2^{t+1+\b}(q_j).
\ees
This proves the statement of the lemma for $\varepsilon=1$ and $\sigma=0$. $\Box$

With each element of $B$ of the form
\bee\label{f7}
u= \delta_1^a(x_1)\delta_2^b(x_2)\delta_1^s\delta_2^t(q_i), \ \ a,b\geq 1, s,t\geq 0,
\eee
we associate the configuration $[i,s,t]$ of the Minsky machine $M$.

Denote by $V_{\varepsilon \sigma}$ the set of all elements of $B^e$ of the form
\bes
w=(\delta_1^a(x_1))^{1-\varepsilon}(\delta_2^b(x_2))^{1-\sigma}\delta_1^s\delta_2^t,
\ees
where $a,b\geq 1$ and $s,t\geq 0$

Every $v\in V$ can be uniquely represented as $v=v_1v_2$, where $v_1$ is a monomial in the variables (\ref{f3}) and $v_2$ is a monomial in the variables (\ref{f5}). We have $\mathrm{Deg}(v)=\mathrm{Deg}(v_1)$ and $\mathrm{Deg}(v_2)=(0,0)$. We denote $v_1$ by $\{v\}$.

\begin{lm}\label{l4}
Let $u$ and $v$ be two elements of the form (\ref{f7}).  Then
\bee\label{f8}
u-v=\overline{w g(i,\varepsilon,\sigma)}
\eee
for some $w\in V_{\varepsilon \sigma}$ if and only if $\{u\}=\{v\}$ and  $[u]\rightarrow [v]$ as a result of the execution of the command (\ref{f1}) (or $[v]\rightarrow [u]$ if $1+1=0$ in $\Phi$).
\end{lm}
\Proof We consider only the case $\varepsilon=1$ and $\sigma=0$. Then
$w\in V_{1 0}$ has the form
\bes
w=\delta_2^r(x_2)\delta_1^s\delta_2^t.
\ees
By Lemma \ref{l3},
\bes
\overline{w(g(i,1,0))}=\delta_2^r(x_2)\overline{\delta_1^s\delta_2^t(g(i,1,0))}\\
=\delta_1^s(x_1)\delta_2^r(x_2)
\delta_2^{t+1}(q_i)
-\delta_1^s(x_1)\delta_2^r(x_2)\delta_1^{\a}
\delta_2^{t+1+\b}(q_j)
\ees

Assume that $1+1\neq 0$ in $\Phi$. Then (\ref{f8}) holds if and only if
\bes
u=\delta_1^s(x_1)\delta_2^r(x_2)
\delta_2^{t+1}(q_i),   \ \ v=\delta_1^s(x_1)\delta_2^r(x_2)\delta_1^{\a}
\delta_2^{t+1+\b}(q_j).
\ees
Notice that $u,v$ has the form (\ref{f7}), $\{u\}=\{v\}$, and $[u]=[i,0,t+1]$ and $[v]=[j,\a,t+1+\b]$. We get $[u]\rightarrow [v]$ as a result of the execution of the command (\ref{f1}).

If $1+1= 0$ in $\Phi$, then
\bes
v=\delta_1^s(x_1)\delta_2^r(x_2)
\delta_2^{t+1}(q_i),   \ \ u=\delta_1^s(x_1)\delta_2^r(x_2)\delta_1^{\a}
\delta_2^{t+1+\b}(q_j)
\ees
is possible. In this case we get $[v]\rightarrow [u]$.
$\Box$

For each $\varepsilon, \sigma= 0,1$, denote by $W_{\varepsilon \sigma}$ the set of all elements of the form
\bes
x_1^{1-\varepsilon}x_2^{1-\sigma}\delta_1^i\delta_2^j, \ \ \ i,j\geq 0,
\ees
such that $i=0$ if $\varepsilon=1$ and $j=0$ if $\sigma=1$. In particular, we have $W_{1 1}=\{1\}$.

\begin{co}\label{c1}
Let $u$ and $v$ be two elements of the form (\ref{f7}) such that $\{u\}=\{v\}=x_1x_2$. Then the equality (\ref{f8}) holds only if $w\in W_{\varepsilon \sigma}$ and in this case
\bes
u-v=w g(i,\varepsilon,\sigma).
\ees
\end{co}
\Proof Again consider only the case $\varepsilon=1$ and $\sigma=0$. If $\{u\}=\{v\}=x_1x_2$, then using the proof of Lemma \ref{l4}, we get
\bes
s=0, r=0,
\ees
and consequently,
\bes
u=x_1x_2
\delta_2^{t+1}(q_i), v=x_1x_2\delta_1^{\a}
\delta_2^{t+1+\b}(q_j), w=x_2\delta_2^t\in W_{1 0}.
\ees
Then
\bes
w(g(i,1,0))=u-v. \ \ \Box
\ees

\begin{co}\label{c2}
If $m\in S$, then $g_m\in \widetilde{I}$.
\end{co}
\Proof
If $m\in S$, then there exists a sequence of configurations
\bes
[1,2^{2^m},0]=c_0\rightarrow c_1\rightarrow\ldots c_r=[0,1,0].
\ees
of the Minsky machine $M$.  For each configuration $c_i$ there exists a unique element $u_i$ of the form (\ref{f7}) such that
$[u_i]=c_i$  and $\{u_i\}=x_1x_2$. Notice that $g_m=u_0-u_r$. By Lemma \ref{l3} and Corollary \ref{c1}, we have $u_i- u_{i+1}\in \widetilde{I}$ for all $0\leq i<r$. Consequently,
\bes
g_m=u_0-u_r=(u_0-u_1)+(u_1-u_2)+\ldots+(u_{r-1}-u_r)\in \widetilde{I}. \ \ \ \Box
\ees

\begin{lm}\label{l5}
If $g_m$ is a linear combination of elements of the form
\bes
w g(i,\varepsilon,\sigma),
\ees
where $w\in W_{\varepsilon \sigma}$, then $m\in S$.
\end{lm}
\Proof Put $u=x_1x_2\delta_1^{2^{2^m}}(q_1)$ and $v=x_1x_2\delta_1(q_0)$. Then $g_m=u-v$. By Lemma \ref{l4}, we may assume that
\bee\label{f9}
u-v=\lambda_1(u_1-v_1)+\lambda_2(u_2-v_2)+\ldots+\lambda_r(u_r-v_r),
\eee
where all $u_i,v_i$ are elements of the form (\ref{f7}) and $[u_i]\rightarrow [v_i]$ for all $1\leq i\leq r$. Assume that $r$ is the minimal number satisfying (\ref{f9}). This condition immediately implies that $u_i\neq v_i$.

In order to prove that $m\in S$, it is sufficient to show the existence of a sequence of configurations of the form
\bes
[u]\rightarrow\ldots\rightarrow [v].
\ees

If $u_i=u_j$, then $v_i=v_j$ since $[u_i]\rightarrow [v_i]$ for all $i$. Consequently, we may assume that all $u_1,u_2,\ldots,u_r$ are different. The machine $M$ at the configuration $[v]=[0,1,0]$ immediately stops its work since it is in the internal state $q_0$. For the same reason, the machine at the configurations $[u_1],[u_2],\ldots,[u_r]$ is not in the internal state $q_0$.
This means that $v$ contains $q_0$ and $u_1,u_2,\ldots,u_r$ do not contain it.

All elements $u,v,u_1,v_1,\ldots,u_r,v_r$ belong to a linearly independent set of elements (\ref{f7}). Then the equality (\ref{f9}) implies that
$v$ coincides with one of $v_1,v_2,\ldots,v_r$. Without loss of generality, we may assume that $v=v_1$. If $u=u_1$, then $[u]=[u_1]\rightarrow [v_1]=[v]$. Otherwise (\ref{f9}) implies that
$u_1$ coincides with one of $v_2,\ldots,v_r$. Without loss of generality, we may assume that $u_1=v_2$.

Continuing this discussion, we may assume
 that $v=v_1, u_1=v_2,\ldots,u_s=v_{s+1}$, $u\neq u_1,u_2,\ldots,u_s$, and $s$ is the maximal number with this property. If $u=u_{s+1}$, then
\bes
[u]=[u_{s+1}]\rightarrow [v_{s+1}]=[u_s]\rightarrow\ldots \rightarrow [v_1]=[v].
\ees
If $u\neq u_{s+1}$, then (\ref{f9}) implies that
$u_{s+1}$ coincides with one of $v_1,v_2,\ldots,v_r$.
If $u_{s+1}=v_j$ for some $1\leq j\leq s+1$, then we get
\bes
(u_j-v_j)+(u_{j+1}-v_{j+1})+\ldots+(u_{s+1}-v_{s+1})=0.
\ees
This allows to reduce the number $r$ in (\ref{f9}).

Consequently, $u_{s+1}$ coincides with one of $v_{s+2},\ldots,v_r$.
We may assume that $u_{s+1}=v_{s+2}$ and this contradicts the maximality of $s$. $\Box$

We intentionally avoided to use the acyclicity of $M$ in the proof of Lemma \ref{l6}. The next lemma is not true for machines without cycles.

\begin{lm}\label{l6} The set of elements of the form
\bee\label{f10}
\overline{w g(i,\varepsilon,\sigma)},
\eee
where $w\in V_{\varepsilon \sigma}$, $1\leq i\leq n$, and $\varepsilon,\sigma\in \{0,1\}$,
are linearly independent over $\Phi$.
\end{lm}
\Proof By Lemma \ref{l4}, every element of the form (\ref{f10}) can be represented as $u-v$, where $u$ and $v$ are elements of the form (\ref{f7}) such that $\{u\}=\{v\}$ and $[u]\rightarrow [v]$.

First of all, notice that $u-v\neq 0$. If $u=v$, then $[u]\rightarrow [v]=[u]$ becomes a nontrivial cycle of the machine $M$. Recall that $M$ is acyclic.

A nontrivial linear dependence of elements of the form
(\ref{f10}) can be written in the form
\bee\label{f11}
\lambda_1(u_1-v_1)+\lambda_2(u_2-v_2)+\ldots+\lambda_r(u_r-v_r)=0,
\eee
where $0\neq \lambda_1,\lambda_2,\ldots,\lambda_r\in \Phi$ and $u_i$ and $v_i$ are elements of the form (\ref{f7}) such that $\{u_i\}=\{v_i\}$ and $[u_i]\rightarrow [v_i]$ for all $i$.

We noticed that $u_i=u_j$ implies $v_i=v_j$. Therefore, we may assume that all $u_1,u_2,\ldots,u_r$ are different. Start with $v_1$ as in the proof of Lemma \ref{l5}. We have $[u_1]\rightarrow [v_1]$.
Then (\ref{f11}) implies that $u_1$ coincides with one of $v_2,v_3,\ldots,v_r$. We may assume that $u_1=v_2$. If $u_2=v_1$, then we get a cycle $[u_2]\rightarrow [v_2]=[u_1]\rightarrow [v_1]=[u_2]$ of the machine $M$. We also know that $u_2\neq v_2$. Consequently, $u_2$ coincides with one of $v_3,\ldots,v_r$.

Assume
 that $u_1=v_2,\ldots,u_s=v_{s+1}$ and $s$ is the maximal number with this property. If $u_{s+1}=v_j$ for some $1\leq j\leq s+1$, the we get a cycle
\bes
[u_{s+1}]\rightarrow [v_{s+1}]=[u_s]\rightarrow[v_s]\rightarrow\ldots \rightarrow [v_j]=[u_{s+1}]
\ees
of $M$. Consequently, $u_{s+1}\neq v_1,v_2,\ldots,v_{s+1}$. Then (\ref{f11}) implies that
$u_{s+1}$ coincides with one of $v_{s+2},\ldots,v_r$. We may assume that $u_{s+1}=v_{s+2}$ and this contradicts to the maximality of $s$. $\Box$

\begin{lm}\label{l7}
If $g_m\in \widetilde{I}$, then $m\in S$.
\end{lm}
\Proof Let $g_m\in \widetilde{I}$. Then
\bee\label{f12}
g_m=\sum_{s,i,\varepsilon,\sigma} \lambda_{s,i,\varepsilon,\sigma} w_s(\varepsilon,\sigma) g(i,\varepsilon,\sigma),
\eee
where $\lambda_{s,i,\varepsilon,\sigma}\in \Phi$ and $w_s(\varepsilon,\sigma)\in B^e$ are elements of the form (\ref{f6}).

Notice that all elements $w_s(\varepsilon,\sigma) g(i,\varepsilon,\sigma)$ and $g_m$ are homogeneous with respect to each set of variables
(\ref{f3})-(\ref{f5}) and with respect to degree function $\deg$. Recall that
\bes
\deg(g_m)=3, \deg(g(i,\varepsilon,\sigma))=1+\varepsilon +\sigma.
\ees
Consequently, we may assume that $w_s(\varepsilon,\sigma)\in V_{\varepsilon,\sigma}$ in (\ref{f12}).

Suppose that there exists at least one $w_s(\varepsilon,\sigma)$ that does not belong to
$W_{\varepsilon,\sigma}$. In this case we get $\mathrm{Deg}(w_s(\varepsilon,\sigma) g(i,\varepsilon,\sigma))>(1,1)$ by Corollary \ref{c1}. Let $(c,d)$ be the highest degree of all elements $w_s(\varepsilon,\sigma) g(i,\varepsilon,\sigma)$ participating in (\ref{f12}) with respect to $\mathrm{Deg}$. We have $(c,d)>(1,1)=\mathrm{Deg}(g_m)$. Then the equality (\ref{f12}) implies a nontrivial linear dependence of the highest homogeneous parts  $\overline{w_s(\varepsilon,\sigma) g(i,\varepsilon,\sigma)}$ of elements  $w_s(\varepsilon,\sigma) g(i,\varepsilon,\sigma)$ with the degree $(c,d)$. It is impossible by Lemma \ref{l6}.

Therefore, we may assume that all $w_s(\varepsilon,\sigma)$ in (\ref{f12}) belong to $W_{\varepsilon,\sigma}$. Then Lemma \ref{l5} implies that $m\in S$. $\Box$

{\bf Proof of Proposition 1}. Notice that $f_m\in I+J$ in the algebra $A$ if and only if $g_m\in \widetilde{I}$ in the algebra $B$. By Corollary \ref{c2} and Lemma \ref{l7}, $g_m\in \widetilde{I}$ if and only if $m\in S$. $\Box$

Proposition \ref{p1} immediately implies the next result.
\begin{theor}\label{t1}
The ideal membership problem for differential polynomial algebras with at least two basic derivations is algorithmically undecidable.
\end{theor}
\Proof Let $S$ be a recursively enumerable but nonrecursive set \cite{Malcev}. This means that there is no algorithm which determines for every natural $m$ whether $m\in S$ or not. Assume that the algebra $A$ and its ideal $I+J$ are constructed by the commands of an acyclic machine $M$ calculating the characteristic function of $S$. By Proposition \ref{p1}, $m\in S$ if and only if $f_m\in I+J$. Consequently, there is no algorithm which determines for all $m$ whether $f_m\in I+J$ or not. $\Box$

\section{The subalgebra membership problem}

\hspace*{\parindent}

Let $\Delta = \{\delta_1,\ldots, \delta_m\}$ be a basic set of derivation operators and $\Theta$ is the free commutative monoid on $\Delta$.
Let $A=\Phi\{x_1,x_2,\ldots,x_n\}$ be a differential polynomial algebra over $\Phi$ and let
$I=[f_1,f_2,\ldots,f_r]$ be a differential ideal of $A$ generated by $f_1,f_2,\ldots,f_r$.
Let $B=\Phi\{x_1,x_2,\ldots,x_n,t\}$ be a differential polynomial algebra with one more free differential variable $t$. We identify $A$ with the corresponding subalgebra of $B$. Denote by $S_I$ the differential subalgebra of $B$ generated by
\bes
x_1,x_2,\ldots,x_n,\delta_1(t),\delta_2(t),\ldots,\delta_m(t),tf_1,tf_2,\ldots,tf_r.
\ees
\begin{pr}\label{p2}
Let $f\in A$. Then $f\in I$ if and only if $tf\in S_I$.
\end{pr}
\Proof Let $M$ be the free commutative monoid generated by all elements
$x_i^{\theta}$, where $1\leq i\leq n$ and $\theta\in \Theta$. Then every element of $A^e$ is a linear combination of elements of the form
\bes
m\theta, \ \ m\in M, \theta\in \Theta.
\ees
If $f\in I$, then
\bes
f=\sum_{m,\theta,i} \lambda_{m,\theta,i}m\theta f_i
\ees
for some $\lambda_{m,\theta,i}\in \Phi$, $m\in M$, $\theta\in \Theta$, and $1\leq i\leq r$.

Let $T$ be the subalgebra of $S_I$ generated by the elements
\bes
x_1,x_2,\ldots,x_n,\delta_1(t),\delta_2(t),\ldots,\delta_m(t).
\ees
Notice that $M\subset A\subset T$. It is easy to check that
\bee\label{f13}
\theta(tf_i)=t\theta(f_i)+g, g\in T.
\eee
Consequently,
\bes
tf=\sum_{m,\theta,i} \lambda_{m,\theta,i}mt\theta(f_i)=\sum_{m,\theta,i} \lambda_{m,\theta,i}m\theta (tf_i)+g, \ \ g\in T.
\ees
Therefore, $tf\in S_I$.

Denote by $\deg_t$ the polynomial degree function on $B$ such that $\deg_t(t^\theta)=1$ and $\deg_t(x_i^\theta)=0$ for all $i$ and $\theta$. All generates of the subalgebra $S_I$ are homogeneous with respect to $\deg_t$. Denote by $H$ the subspace of all homogeneous elements of degree $1$ of $S_I$ with respect to $\deg_t$. Then every element of $H$ is a linear combination of elements of the form
\bes
mt^{\theta_1}, m \theta(tf_i),
\ees
where $\theta_1,\theta\in \Theta$, $|\theta_1|\geq 1$, $m\in M$, and $1\leq i\leq r$.
Taking into account (\ref{f13}), we may assume that every element of $H$ is a linear combination of elements
\bes
mt^{\theta_1}, mt \theta(f_i).
\ees
Moreover, every element of $H$ divisible by $t$ is a linear combination of elements
\bes
mt \theta(f_i), \ \ m\in M, \theta\in \Theta, 1\leq i\leq r.
\ees

Assume that $tf\in S_I$ for some $f\in A$. We have $tf\in H$ and $tf$ divisible by $t$. Then
\bes
tf=\sum_{m,\theta,i} \lambda_{m,\theta,i} mt \theta(f_i).
\ees
Consequently,
\bes
f=\sum_{m,\theta,i} \lambda_{m,\theta,i} m \theta(f_i)\in I. \ \ \Box
\ees

\begin{theor}\label{t2}
The subalgebra membership problem for differential polynomial algebras with at least two basic derivations is algorithmically undecidable.
\end{theor}
\Proof By Theorem \ref{t1}, the ideal membership problem is undecidable. By Proposition \ref{p2}, if the ideal membership problem for $A$ is undecidable, then the subalgebra membership problem for $B$ is also undecidable. $\Box$

The method of this section gives the undecidability of the subalgebra membership problem for free algebras of many varieties of algebras \cite{U95A}. But this method does not work for subfields of fields. Recall that the subfield membership problem for fields of rational functions $k(x_1,x_2,\ldots,x_n)$ over a constructive field $k$ is also positively solved by means of Gr\"obner bases \cite{Sweedler93}. The subfield membership problem for differential fields of rational functions remains open.

\end{document}